\documentclass[12pt]{article}
\usepackage{amsxtra,amssymb,amsthm,amsmath,latexsym}

\theoremstyle{plain}

\newtheorem{theorem}{Theorem}[section]

\newtheorem*{remark3.2}{Remark 3.2}
\newtheorem*{remark3.3}{Remark 3.3}
\numberwithin{equation}{section}

\def\nd{\noindent}

\def\oH{\buildrel\circ\over H}
\def\oH1{\buildrel\circ\over H\kern-.02in{}^1}

\begin{document}


\title{                  
On a new notion of regularizer
  \thanks{key words: ill-posed problems, regularizer, stable numerical
differentiation, new notion of regularizer
    }
}

\author{
A.G. Ramm\\
LMA/CNRS, 31 Chemin Joseph Aiguier, Marseille 13402, France\\
and Mathematics Department, Kansas State University, \\
 Manhattan, KS 66506-2602, USA\\
ramm@math.ksu.edu\\
}

\date{}

\maketitle\thispagestyle{empty}

\begin{abstract}

A new understanding of the notion of regularizer is proposed. It is
argued that this new notion is more realistic than the old one
and better fits the practical computational needs. An example
of the regularizer in the new sense is given. A method for
constructing regularizers in the new sense is proposed and justified.

\end{abstract}


\section{Introduction}
Let
\begin{equation} \label{1.1} A(u)=g, \end{equation}
where $A:X\to Y$ is a closed,
possibly nonlinear, map from a Banach space $X$ into a Banach space $Y$.
Problem \eqref{1.1} is called ill-posed if $A$ is not a homeomorphism
of $X$ onto $Y$, that is, either equation \eqref{1.1} does not have a solution,
or the solution is non-unique, or the solution does not depend on $g$
continuously. Let us assume that \eqref{1.1} has a solution $u$ and this
solution is unique but $A^{-1}$ is not continuous.
Given noisy data $g_\delta$, $\|g_\delta-g\|\leq \delta$,
one wants to construct a stable approximation $u_{\delta}$ of the solution
$u$, $\|u_\delta-u\|\to 0$ as $\delta \to 0$.
This is often done with the help of a regularizer.
Traditionally (see, e.g., \cite{2}) one calls a family of
operators $R_h$ a regularizer for problem \eqref{1.1} if:

\nd a) $R_hA(u)\to u$ as $ h\to \infty$ for any $u\in D(A)$,\\
\nd b) $R_hg_\delta$ is defined for any $g_\delta\in Y$
and there exists $h(\delta)\to 0$ as $\delta\to 0$ such that
\begin{equation} \label{ast} \|R_{h(\delta)} g_\delta-u\| \to 0 
\tag{{$\ast$}} \hbox { as } \delta\to 0, \end{equation}
where $u$ solves \eqref{1.1}.

In this definition $u$ is fixed and \eqref{ast} must hold for any
$g_\delta\in B(g,\delta):=\left\{ g_\delta:\|g_\delta-g\|\leq \delta 
\right\}$.

In practice one does not know the solution $u$.
The only available information is a family $g_\delta$ and some
priori information about the solution $u$.
This information very often
consists of the knowledge that $u\in {\mathcal K}$, 
where ${\mathcal K}$ is a compactum in $X$.
Thus
$u\in S_\delta
  :=\left\{ v:\|A(v)-g_\delta\|\leq\delta,\ v\in {\mathcal K}\right\}$.
We assume that the operator $A$ is known exactly, and we always
assume that $g_\delta\in B(g,\delta)$, where $g=A(u)$.

It is natural to call a family of operators $R(\delta)$ a regularizer
(in the {\bf new sense}) if
\begin{equation} \label{1.2}
 \sup_{v\in S_\delta} \|R(\delta) g_\delta-v\|
 \leq \eta(\delta)\to 0 \hbox{\ as\ } \delta \to 0. \end{equation}
There is a crucial difference between our definition \eqref{1.2}
and the standard definition \eqref{ast}:
in \eqref{ast} $u$ is fixed,
while in \eqref{1.2} $v$ is any element of $S_\delta$ and the supremum, 
over all such $v$, of the norm in \eqref{1.2} must go to zero as 
$\delta\to 0$.

The new definition is more realistic and fits more computational needs
because not only the solution $u$ to \eqref{1.1} satisfies the inequality
$\|A(u)-g_\delta\|\leq \delta$,
but many $v$ satisfy such an inequality
$\|Av-g_\delta\|\leq \delta$, $v\in {\mathcal K}$,
and 
the data $g_\delta$ may correspond,
in fact, to any $v\in S_\delta$, and not only to the solution
of equation (1.1).
Therefore it is more natural  to use definition \eqref{1.2} than \eqref{ast}.

Our aim is to illustrate the practical difference in these two
definitions by an example, and to construct regularizer in the sense
\eqref{1.2} for problem \eqref{1.1} under the following assumptions:

\nd A1) $A:X\to Y$ is a closed, possibly nonlinear, injective map,
$g\in {\mathcal R}(A)$, 

and

\nd A2) $\phi: D(\phi)\to [0,\infty)$, $\phi(u)>0$ if $u\not= 0$,
domain $D(\phi)\subseteq D(A)$,
the set ${\mathcal K}={\mathcal K}_c:=\{v:\phi(v)\leq c\}$
is compact in $X$ and contains a sequence
$v_n\to v$,
and
$\phi(v)\leq \liminf_{n\to\infty} \phi(v_n)$.

The last inequality holds if $\phi$ is lower semicontinuous.
Any Hilbert space norm and norms in reflexive Banach spaces
have this property.

Examples in which assumptions A1) and A2) are satisfied are numerous.

\nd{\bf Example 1.}
$A$ is a linear injective compact operator, $g\in{\mathcal R}(A)$,
$\phi(v)$ is a norm on $X_1\subset X$, where $X_1$ is densely
imbedded in $X$, the embedding $i:X_1\to X$ is compact, and
$\phi(v)$ is lower semicontinuous.

\nd{\bf Example 2.}
$A$ is a nonlinear injective continuous operator $g\in{\mathcal R}(A)$,
$A^{-1}$ is not continuous, $\phi$ is as in Example 1.

\nd{\bf Example 3.}
$A$ is linear, injective, densely defined, closed operator,
$g\in{\mathcal R}(A)$, $A^{-1}$ is unbounded, $\phi$ is as in Example 1,
$X_1\subseteq D(A)$.

In section 2 
it is shown that a regularizer in
the sense \eqref{ast} may be not a regularizer in the sense \eqref{1.2}.
In section 3 a theoretical construction of a regularizer in the sense
\eqref{1.2} is given.

\section{Example: stable numerical differentiation.}
Here we use the results from \cite{R1} - \cite{R7}.

Consider stable numerical differentiation of noisy data.
The problem is:
\begin{equation} \label{2.1}
Au:=\int^x_0 u(s)\,ds=g(x), \quad g(0)=0, \quad 0\leq x\leq 1.
\end{equation}
The data are: $g_\delta$ and $M_a$, where $\|g_\delta-g\|\leq \delta$,
the norm is $L^\infty(0,1)$ norm, and
$\|u\|_a\leq M_a$, $a\geq 0$.
The norm
$$\|u\|_a := \sup_{\substack {x,y\in[0,1] \\ x\neq y}}
  \ \frac{|u(x)-u(y)|}{|x-y|^a}
 + \sup_{0\leq x\leq 1} |u(x)|\quad \hbox{\ if\ }  0\leq a\leq 1,
$$
$$\|u\|_a:=\sup_{0\leq x\leq 1} (|u(x)| + |u'(x)|)
 +\sup_{\substack {x,y\in[0,1] \\ x\neq y}}
  \ \frac{|u'(x)-u'(y)|}{|x-y|^{a-1}},
 \quad 1<a\leq 2.
$$
If $a>1$, then we define
\begin{equation} \label{2.2}
R(\delta)g_\delta:=
\begin{cases}
\frac{g_\delta(x+h(\delta))-g_\delta(x-h(\delta))}{2h(\delta)},
   & h(\delta)\leq x\leq 1-h(\delta),\\
\frac{g_\delta(x+h(\delta))-g_\delta(x)}{h(\delta)},
   & 0\leq x< h(\delta), \\
\frac{g_\delta(x)-g_\delta(x-h(\delta))}{h(\delta)},
   & 1-h(\delta)<x\leq 1,
\end{cases} \end{equation}
where
\begin{equation}
\label{2.3} h(\delta)=c_a\delta^{\frac{1}{a}},
\end{equation}
and $c_a$ is a constant given explicitly (cf \cite{R2}). 

We prove that
\eqref{2.2} is a regularizer for \eqref{2.1} in the sense \eqref{1.2},
and ${\mathcal K}:=\{v:\|v\|_a \leq M_a,\ a>1\}$.
In this example we do not use lower semicontinuity of the norm $\phi(v)$
and do not define $\phi$. 

Let $S_{\delta,a}:=\{ v:\|Av-g_\delta\| \leq \delta, \ \|v\|_a \leq M_a 
\}$. 
To prove that \eqref{2.2}-\eqref{2.3} is a regularizer in the sense
\eqref{1.2} we use the estimate
\begin{equation} \label{2.4}
\begin{aligned}
\sup_{\substack {v\in S_{\delta,a}}}  &  \|R(\delta) g_\delta-v\|
\leq \sup_{\substack {v\in S_{\delta,a}}}
  \{ \|R(\delta)(g_\delta-Av)\| +\| R(\delta)Av-v\| \}
\leq \frac{\delta}{h(\delta)} +M_a h^{a-1}(\delta) \leq
 \\&
\leq c_a \delta^{1-\frac{1}{a}} :=\eta(\delta)\to 0
  \hbox{\ as\ } \delta \to 0. 
\end{aligned}\end{equation}
Thus we have proved that \eqref{2.2}-\eqref{2.3} is a regularizer
in the sense \eqref{1.2}.

If $a=1$, and $M_1<\infty$, then one can prove that {\bf there is no}
regularizer for problem \eqref{2.1} in the sense \eqref{1.2}
even if the regularizer is sought in the set of all operators,
including nonlinear ones. More precisely, it is proved in \cite{R3}
(see also  \cite{R5})
that
$$\inf_{R(\delta)} \ \sup_{v\in S_{\delta,1}}
 \ \|R(\delta) g_\delta - v\| \geq c>0,$$
where $c>0$ is a constant independent of $\delta$
and the infimum is taken over all operators
$R(\delta)$ acting from $L^{\infty}(0,1)$
into  $L^{\infty}(0,1)$, including nonlinear ones.
On the other hand, if $a=1$ and $M_1<\infty$, then
a regularizer in the sense \eqref{ast} does exist,
but the rate of convergence in (*) may be as slow
as one wishes, if
$u(x)$ is chosen suitably \cite{R6}, \cite{R7}.

\section{Construction of a regularizer in the sense \eqref{1.2}}
Assuming A1) and A2) (see Section 1) throughout this Section, let us 
construct a regularizer
for \eqref{1.1} in the sense \eqref{1.2}. We use the ideas from
\cite{R7}.
Define $F_\delta(v):=\|Av-g_\delta\| + \delta\phi(v)$ and
consider the minimization problem of finding the infimum $m(\delta)$ of 
the 
functional $ F_\delta(v)$ on a set $ S_\delta$, defined below:
\begin{equation} \label{3.1}
 F_\delta(v)=\inf:=m(\delta),
 \ v\in S_\delta:=\{v:\|Av-g_\delta\| \leq \delta,
 \ \phi(v)\leq c\},
\end{equation}
so that ${\mathcal K}={\mathcal K}_c:=\{ v:\phi(v)\leq c\}$.
The constant $c>0$ can be chosen arbitrary large and fixed at
the beginning of the argument, and then one can choose a smaller
constant $c_1$, specified below.
Since $F_\delta(u)=\delta + \delta\phi(u):=c_1\delta$,
$c_1:=1+\phi(u)$, where $u$ solves \eqref{1.1}, one concludes that
\begin{equation} \label{3.2}
 m(\delta)\leq c_1\delta.
\end{equation}
Let $v_j$ be a minimizing sequence, such that
$F_\delta(v_j)\leq 2m(\delta)$.
Then $\phi(v_j)\leq 2 c_1$. By assumption A2),
as $j\to\infty$, one has:
\begin{equation}  \label{3.3}
 v_j\to v_\delta,
 \, \phi(v_\delta) \leq 2c_1.
\end{equation}
Take $\delta=\delta_m\to 0$ and denote
$v_{\delta_m}:=w_m$.
Then \eqref{3.3} and assumption A2) imply existence
of a subsequence, denoted again $w_m$, such that:
\begin{equation} \label{3.4}
 w_m\to w,\ A(w_m)\to A(w), \ \|A(w)-g\|=0.
\end{equation}
Thus $A(w)=g$ and, since $A$ is injective, $w=u$, where $u$ is the
unique solution to \eqref{1.1}.

Define now $R(\delta)g_\delta$ by the formula 
$R(\delta)g_\delta:=v_\delta$.

\begin{theorem}
$R(\delta)$ is a regularizer for problem \eqref{1.1} in the
sense \eqref{1.2}
\end{theorem}

\begin{proof}
Assume the contrary:
\begin{equation}\label{3.5}
  \sup_{v\in S_\delta} \|R(\delta)g_\delta - v\|
  = \sup_{v\in S_\delta} \|v_\delta-v\| \geq \gamma >0,
\end{equation}
where $\gamma>0$ is a constant independent of $\delta$.
Since $\phi(v_\delta)\leq 2c_1$ by \eqref{3.3}, and
$\phi(v)\leq c$, one can choose convergent in $X$ sequences
$w_m:=v_{\delta_m}\to \tilde w$, $\delta_m\to 0$, and $v_m\to\tilde 
v$,
such that $\|w_m-v_m\|\geq \frac{\gamma}{2}$,
$\| \tilde w -\tilde v\|\geq\frac{\gamma}{2}$,
and $A(\tilde w)=g$, $A(\tilde v)=g$.
Therefore, by the injectivity of $A$, $\tilde w=\tilde v=u,$ and one gets 
a contradiction
with the inequality $\|\tilde w-\tilde v\|\geq\frac{\gamma}{2}>0$.
This contradiction proves Theorem 3.1.    

Note that the conclusions $A(\tilde w)=g$ and $A(\tilde v)=g$ follow
from the inequalities $\|A(v_\delta)-g_\delta\|\leq \delta$
and $\|A(v)-g_\delta\| \leq \delta$
after passing to the limit $\delta\to 0$, using assumption A2).
\end{proof}

\begin{remark3.2}
Our argument is closely related to a generalization
of a classical result \cite{DS} which says that an injective and
continuous map of a compactum in a Banach space  has continuous
inverse. The generalization of this result, given in \cite{R4} (Chapter 5, 
Section 6, Lemma 2),
says that the same conclusion holds if $A$ is an injective
and closed map. This generalization, in other words, 
claims that if $A$ is an injective and closed map of a
compactum ${\mathcal K}$ in a Banach space
into a Banach space, then the
modulus of continuity $\omega(\delta)$ of the inverse operator
$A^{-1}$ on the set $A({\mathcal K})$ tends to zero:
$$
\sup_{\substack{\|Av-Aw\|\leq \delta \\ v,w \in{\mathcal K}}}
\|v-w\|:=\omega(\delta)\to 0 \hbox { as } \delta \to 0. 
$$
Therefore, in our proof of Theorem 3.1, one has
$\sup_{v\in S_\delta}\|v_\delta-v\|\to 0$ as $\delta\to 0$.

Indeed, $v\in S_\delta$ and $v_\delta\in S_\delta$,
so $\|A(v)-g_\delta\|\leq \delta$ and $\|A(v_{\delta})-g_\delta\|\leq 
\delta$. Thus $\|A(v)-A(v_\delta)\|\leq 2\delta.$
Since $v_\delta$ and $v$ belong to a compactum
${\mathcal K}_{c_1}:=\{ v:\phi(v)\leq 2c_1\}$,
Theorem 3.1 follows from the property $\omega(\delta)\to 0$
as $\delta\to 0$.
\end{remark3.2}

\begin{remark3.3}
We have assumed that $A$ is injective, that is, equation \eqref{1.1}
has at most one solution. This assumption can be relaxed: one may
assume that  $u$ solves \eqref{1.1} and there is an $r>0$,
such that in the ball $B(u,r):=\{v:\|v-u\| \leq r\}$ the solution $u$ is
unique, but globally \eqref{1.1} may have many solutions.
Our arguments remain valid if the compactum ${\mathcal K}:=
\{v:\phi(v)\leq c,\ v\in B(u,r)\}$.

One may drop the injectivity of $A$ assumption in A1) in Section 1.
In this case (3.4) yields $w\in U$, where $U:=\{w: A(w)=g\}$. Thus,
$\lim_{\delta \to 0}\rho( R(\delta)g_\delta, U)=0$, where
$\rho(w,U)$ is the distance from an element $w$ to the set
$U$.

In applications, when a physical problem is reduced to equation (1.1),
one wants to have a unique solution to (1.1). If the solution is 
non-unique, that is, $U$ contains more than one element, then one
wants to impose additional conditions which select a unique solution in 
the set $U$, thus making the operator $A$ injective.
\end{remark3.3}

\end{document}